# LIKELIHOOD BASED INFERENCE FOR MONOTONE RESPONSE MODELS[1]

By Moulinath Banerjee

*University of Michigan*

The behavior of maximum likelihood estimates (MLEs) and the likelihood ratio statistic in a family of problems involving pointwise nonparametric estimation of a monotone function is studied. This class of problems differs radically from the usual parametric or semi-parametric situations in that the MLE of the monotone function at a point converges to the truth at rate $n^{1/3}$ (slower than the usual $\sqrt{n}$ rate) with a non-Gaussian limit distribution. A framework for likelihood based estimation of monotone functions is developed and limit theorems describing the behavior of the MLEs and the likelihood ratio statistic are established. In particular, the likelihood ratio statistic is found to be asymptotically pivotal with a limit distribution that is no longer $\chi^2$ but can be explicitly characterized in terms of a functional of Brownian motion. Applications of the main results are presented and potential extensions discussed.

**1. Introduction.** A common problem in nonparametric statistics is the need to estimate a function, like a density, a distribution, a hazard or a regression function. Background knowledge about the statistical problem can provide information about certain aspects of the function of interest, which, if incorporated in the analysis, enables one to draw meaningful conclusions from the data. Often, this manifests itself in the nature of shape restrictions (on the function). Monotonicity, in particular, is a shape restriction that shows up very naturally in different areas of application like reliability, renewal theory, epidemiology and biomedical studies. Consequently, monotone functions have been fairly well studied in the literature and several authors have addressed the problem of maximum likelihood estimation under monotonicity constraints. We point out some of the well-known ones.

Received May 2003; revised May 2006.

[1]Supported by NSF Grant DMS-03-06235.

One of the earliest results of this type goes back to Prakasa Rao [21], who derived the asymptotic distribution of the Grenander estimator (the MLE of a decreasing density); Brunk [4] explored the limit distribution of the MLE of a monotone regression function, Groeneboom and Wellner [9] studied the limit distribution of the MLE of the survival time distribution with current status data, Huang and Zhang [14] and Huang and Wellner [13] obtained the asymptotics for the MLE of a monotone density and a monotone hazard with right censored data, and Wellner and Zhang [27] deduced the large sample theory for a pseudo-likelihood estimator for the mean function of a counting process. A common feature of these monotone function problems that sets them apart from the spectrum of regular parametric and semiparametric problems is the slower rate of convergence $(n^{1/3})$ of the maximum likelihood estimates of the value of the monotone function at a fixed point (recall that the usual rate of convergence in regular parametric/semiparametric problems is $\sqrt{n}$). What happens in each case is the following: If $\hat{\psi}_n$ is the MLE of the monotone function $\psi$, then provided that $\psi'(z)$ does not vanish,

$$(1.1) \qquad n^{1/3}(\hat{\psi}_n(z) - \psi(z)) \to_d C(z)\mathbb{Z},$$

where the random variable $\mathbb{Z}$ is a symmetric (about 0) but *non-Gaussian* random variable and $C(z)$ is a constant depending upon the underlying parameters in the problem and the point of interest $z$. In fact, $\mathbb{Z} = \arg\min_h(W(h) + h^2)$, where $W(h)$ is standard two-sided Brownian motion on the line. The distribution of $\mathbb{Z}$ was analytically characterized by Groeneboom [7] and more recently its distribution and functionals thereof have been computed by Groeneboom and Wellner [10].

In this paper we study a class of conditionally parametric models, of the covariate-response type, where the conditional distribution of the response given the covariate comes from a regular parametric model, with the parameter being given by a monotone function of the covariate. We call these *monotone response models*. Here is a formal description: Let $\{p(x,\theta) : \theta \in \Theta\}$, with $\Theta$ being an open subinterval of $\mathbb{R}$, be a one-parameter family of probability densities with respect to a dominating measure $\mu$. Let $\psi$ be an increasing or decreasing continuous function defined on an interval $\tilde{I}$ and taking values in $\Theta$. Consider i.i.d. data $\{(X_i, Z_i)\}_{i=1}^n$ where $Z_i \sim p_Z$, $p_Z$ being a Lebesgue density defined on $\tilde{I}$ and $X_i | Z_i = z \sim p(x, \psi(z))$. Interest focuses on estimating the function $\psi$, since it captures the nature of the dependence between the response $(X)$ and the covariate $(Z)$. If the parametric family of densities, $p(x, \theta)$, is parametrized by its mean, then $\psi(z) = E(X|Z=z)$ is precisely the regression function. In this paper, we study the asymptotics of the MLE of $\psi$ and also the likelihood ratio statistic for testing $\psi$ at a fixed point of interest, with a view to obtaining pointwise confidence sets for $\psi$ of an assigned level of significance. Before we discuss this further, here are some motivating examples to illustrate the above framework.

(a) Consider, for example, the *monotone regression model* where $X_i = \psi(Z_i) + \varepsilon_i$, $\{(\varepsilon_i, Z_i)\}_{i=1}^n$ are i.i.d. random variables, $\varepsilon_i$ is independent of $Z_i$, each $\varepsilon_i$ has mean 0 and variance $\sigma^2$, each $Z_i$ has a Lebesgue density $p_Z(\cdot)$ and $\psi$ is a monotone function. The above model and its variants have been fairly well studied in the literature on isotonic regression (see, e.g., [4, 12, 16, 18]). Now suppose that the $\varepsilon_i$'s are Gaussian. We are then in the above framework: $Z \sim p_Z(\cdot)$ and $X|Z = z \sim N(\psi(z), \sigma^2)$. We want to estimate $\psi$ and test $\psi(z_0) = \theta_0$ for an interior point $z_0$ in the domain of $\psi$.

(b) Another example is the *binary choice model* where we have a dichotomous response variable $X = 1$ or $0$ and a continuous covariate $Z$ with a Lebesgue density $p_Z(\cdot)$ such that $P(X = 1|Z) \equiv \psi(Z)$ is a smooth function of $Z$. In a biomedical context one could think of $X$ as representing the indicator of a disease/infection and $Z$ the level of exposure to a toxin, or the measured level of a biomarker that is predictive of the disease/infection. In such cases it is often natural to impose a monotonicity assumption on $\psi$. A special version of this model is the case 1 interval censoring/current status model that is used extensively in epidemiology and has received much attention among biostatisticians and statisticians (see, e.g., [3, 6, 9, 11]).

(c) The *Poisson regression model* used for count data provides yet another example. Suppose that $Z \sim p_Z(\cdot)$ and $X|Z = z \sim \text{Poisson}(\psi(z))$ where $\psi$ is a monotone function. Here one can think of $Z$ as the distance of a region from a point source (e.g., a nuclear processing plant) and $X$ the number of cases of disease incidence at distance $Z$. The expected number of disease cases at distance $z$ from the source ($\psi(z)$) may be expected to be monotone decreasing in $z$. Variants of this model have received considerable attention in epidemiological contexts [5, 17, 23].

A common feature of all three models described above is the fact that the conditional distribution of the response comes from a one parameter full rank exponential family [in (a), the variance $\sigma^2$ needs to be held fixed]. Our last example below considers a curved exponential family model for the response and is of a fundamentally different flavor in that explicit characterizations of maximum likelihood estimates of $\psi$ are not available in this model, in contrast to the preceding ones.

(d) *Conditional normality under a mean–variance relationship.* Consider the scenario where $Z$ has a Lebesgue density concentrated on an interval $[a, b]$ (with $0 < a < b$) and given $Z = z$, $X \sim p(x, \psi(z))$ for an increasing function $\psi$, with $p(x, \theta)$ being the normal density, $\mu = c\theta^{-2m+1}$ and $\sigma^2 = d\theta^{-2m}$ for some real $m \geq 1$, and $\theta, c, d > 0$. For $m = 1$, this reduces to a normal density with a linear relationship between the mean and the standard deviation. Such a model could be postulated in a real-life setting based on, say, exploratory plots of the mean–variance relationship using observed data, or background knowledge.

Based on existing work, one would expect $\hat{\psi}_n(z_0)$, the MLE of $\psi$ at a prefixed point $z_0$ to satisfy (1.1), with $z$ replaced by $z_0$. As will be seen, this indeed happens. This result permits the construction of (asymptotic) confidence intervals for $\psi(z_0)$ using the quantiles of $\mathbb{Z}$, which are well tabulated. The constant $C(z_0)$ however needs to be estimated and involves nuisance parameters depending on the underlying model, and in particular, the derivative of $\psi$ at $z_0$, estimating which is a tricky affair. Another likelihood-based method of constructing confidence sets for $\psi(z_0)$ would involve testing a null hypothesis of the form $H_{0,\theta} : \psi(z_0) = \theta$, using the likelihood ratio test, for different values of $\theta$, and then inverting the acceptance region of the likelihood ratio test; in other words, the confidence set for $\psi(z_0)$ is formed by compiling all values of $\theta$ for which the likelihood ratio statistic does not exceed a critical threshold. The threshold depends on $0 < \alpha < 1$, where $1 - \alpha$ is the level of confidence being sought, and the asymptotic distribution of the likelihood ratio statistic when the null hypothesis is correct. Thus, we are interested in studying the asymptotics of the likelihood ratio statistic for testing the true (null) hypothesis $H_{0,\theta_0} : \psi(z_0) = \theta_0$. Pointwise null hypotheses of this kind are very important from the *perspective of estimation* since they serve as a conduit for setting confidence limits for the value of $\psi$, through inversion.

A question that arises naturally is whether, similar to the classical parametric case, we can find a *universal limit distribution* for the likelihood ratio statistic when the null hypothesis $H_{0,\theta_0}$ holds, for the monotone response models introduced above. The hope that a universal limit may exist is bolstered by the work of Banerjee and Wellner [3], who studied the limiting behavior of the likelihood ratio statistic for testing the value of the distribution function $(F)$ of the survival time at a fixed point in the current status model. They found that in the limit the likelihood ratio statistic behaves like $\mathbb{D}$, which is a well-defined functional of $W(t) + t^2$ (and is described below). We will show that for our monotone response models, $\mathbb{D}$ does indeed arise as the universal limit law of the likelihood ratio statistic.

We are now in a position to describe the agenda for this paper. In Section 2 we give regularity conditions on the monotone response models under which the results in this paper are developed. We state and prove the main theorems describing the limit distributions of the MLEs and the likelihood ratio statistic. Section 3 discusses applications of the main theorems and Section 4 contains some concluding remarks. The Appendix contains the proofs of some of the lemmas used to establish the main results in Section 2.

**2. Model assumptions, characterizations of estimators and main results.** Consider the general monotone response model introduced in the previous section. Let $z_0$ be an interior point of $\tilde{I}$ at which one seeks to estimate $\psi$.

Assume that: (a) $p_Z$ is positive and continuous in a neighborhood of $z_0$, and (b) $\psi$ is increasing and continuously differentiable in a neighborhood of $z_0$ with $\psi'(z_0) > 0$.

The joint density of the data vector $\{(X_i, Z_i)\}_{i=1}^n$ (with respect to an appropriate dominating measure) can be written as

$$p_n(\psi, \{(X_i, Z_i)\}_{i=1}^n) = \prod_{i=1}^n p(X_i, \psi(Z_i)) \times \prod_{i=1}^n p_Z(Z_i).$$

The second factor on the right-hand side of the above display does not involve $\psi$ and hence is irrelevant as far as computation of MLEs is concerned. Absorbing this into the dominating measure, the likelihood function is given by the first factor on the right-hand side of the display above. Denote by $\hat{\psi}_n$ the unconstrained MLE of $\psi$ and by $\hat{\psi}_n^0$ the MLE of $\psi$ under the constraint imposed by the null hypothesis $H_0: \psi(z_0) = \theta_0$. We assume:

(A.0) With probability increasing to 1 as $n \to \infty$, the MLEs $\hat{\psi}_n$ and $\hat{\psi}_n^0$ exist.

Consider the likelihood ratio statistic for testing the hypothesis $H_0: \psi(z_0) = \theta_0$, where $\theta_0$ is an interior point of $\Theta$. Denoting the likelihood ratio statistic by $2 \log \lambda_n$, we have

$$2 \log \lambda_n = 2 \log \prod_{i=1}^n p(X_i, \hat{\psi}_n(Z_i)) - 2 \log \prod_{i=1}^n p(X_i, \hat{\psi}_n^0(Z_i)).$$

In what follows, assume that the null hypothesis $H_0$ holds.

*Further assumptions.* We now state our assumptions about the parametric model $p(x, \theta)$.

(A.1) The set $\mathcal{X}_\theta = \{x : p(x, \theta) > 0\}$ does not depend on $\theta$ and is denoted by $\mathcal{X}$.
(A.2) $l(x, \theta) = \log p(x, \theta)$ is at least three times differentiable with respect to $\theta$ and is strictly concave in $\theta$ for every fixed $x$ in $\mathcal{X}$. The first, second and third partial derivatives of $l(x, \theta)$ with respect to $\theta$ will be denoted by $\dot{l}(x, \theta), \ddot{l}(x, \theta)$ and $l'''(x, \theta)$.
(A.3) If $T$ is any statistic such that $E_\theta(|T|) < \infty$, then

$$\frac{\partial}{\partial \theta} \int_{\mathcal{X}} T(x) p(x, \theta) \, dx = \int_{\mathcal{X}} T(x) \frac{\partial}{\partial \theta} p(x, \theta) \, dx$$

and

$$\frac{\partial^2}{\partial \theta^2} \int_{\mathcal{X}} T(x) p(x, \theta) \, dx = \int_{\mathcal{X}} T(x) \frac{\partial^2}{\partial \theta^2} p(x, \theta) \, dx.$$

Under these assumptions, $I(\theta) \equiv E_\theta(\dot{l}(X, \theta)^2) = -E_\theta(\ddot{l}(X, \theta))$.

(A.4) $I(\theta)$ is finite and continuous at $\theta_0$.
(A.5) There exists a neighborhood $\mathcal{N}$ of $\theta_0$ such that for all $x$, $\sup_{\theta \in \mathcal{N}} |l'''(x,\theta)| \leq B(x)$ and $\sup_{\theta \in \mathcal{N}} E_\theta(B(X)) < \infty$.
(A.6) The functions $f_1(\theta_1, \theta_2) = E_{\theta_1}(\dot{l}(X,\theta_2)^2)$ and $f_2(\theta_1, \theta_2) = E_{\theta_1}(\ddot{l}(X,\theta_2))$ are continuous in a neighborhood of $(\theta_0, \theta_0)$. Also, the function $f_3(\theta_1, \theta_2) = E_{\theta_1}(\ddot{l}(X,\theta_2)^2)$ is uniformly bounded in a neighborhood of $(\theta_0, \theta_0)$.
(A.7) Set $H(\theta, M)$ to be

$$E_\theta[(|\dot{l}(X,\theta)|^2 + \ddot{l}(X,\theta)^2)(1\{|\dot{l}(X,\theta)| > M\} + 1\{|\ddot{l}(X,\theta)| > M\})].$$

Then $\lim_{M \to \infty} \sup_{\theta \in \mathcal{N}} H(\theta, M) = 0$.

We are interested in describing the asymptotic behavior of the MLEs of $\hat{\psi}_n$ and $\hat{\psi}_n^0$ in local neighborhoods of $z_0$ and that of the likelihood ratio statistic $2 \log \lambda_n$. In order to do so, we first need to introduce the basic spaces and processes (and relevant functionals of the processes) that will figure in the asymptotic theory.

First, define $\mathcal{L}$ to be the space of locally square integrable real-valued functions on $\mathbb{R}$ equipped with the topology of $L^2$ convergence on compact sets. Thus $\mathcal{L}$ comprises all functions $\phi$ that are square integrable on every compact set and $\phi_n$ is said to converge to $\phi$ if $\int_{[-K,K]} (\phi_n(t) - \phi(t))^2 \, dt \to 0$ for every $K$. The space $\mathcal{L} \times \mathcal{L}$ denotes the Cartesian product of two copies of $\mathcal{L}$ with the usual product topology. Also, define $B_{\text{loc}}(\mathbb{R})$ to be the set of all real-valued functions defined on $\mathbb{R}$ that are bounded on every compact set, equipped with the topology of uniform convergence on compacta. Thus $h_n$ converges to $h$ in $B_{\text{loc}}(\mathbb{R})$ if $h_n$ and $h$ are bounded on every compact interval $[-K, K]$ ($K > 0$) and $\sup_{x \in [-K,K]} |h_n(x) - h(x)| \to 0$ for every $K > 0$.

For a real-valued function $f$ defined on $\mathbb{R}$, let $\text{slogcm}(f, I)$ denote the left-hand slope of the GCM (greatest convex minorant) of the restriction of $f$ to the interval $I$. We abbreviate $\text{slogcm}(f, \mathbb{R})$ to $\text{slogcm}(f)$. Also define:

$$\text{slogcm}^0(f) = (\text{slogcm}(f, (-\infty, 0]) \wedge 0) 1_{(-\infty, 0]}$$
$$+ (\text{slogcm}(f, (0, \infty)) \vee 0) 1_{(0, \infty)}.$$

For positive constants $c$ and $d$ define the process $X_{c,d}(z) = cW(z) + dz^2$, where $W(z)$ is standard two-sided Brownian motion starting from 0. Set $g_{c,d} = \text{slogcm}(X_{c,d})$ and $g_{c,d}^0 = \text{slogcm}^0(X_{c,d})$. It is known that $g_{c,d}$ is a piecewise constant increasing function, with finitely many jumps in any compact interval. Also $g_{c,d}^0$, like $g_{c,d}$, is a piecewise constant increasing function, with finitely many jumps in any compact interval and differing, almost surely, from $g_{c,d}$ on a finite interval containing 0. In fact, with probability 1, $g_{c,d}^0$ is identically 0 in some random neighborhood of 0, whereas $g_{c,d}$ is almost surely nonzero in some random neighborhood of 0. Also, the length of the interval $D_{c,d}$ on which $g_{c,d}$ and $g^0_{c,d}$ differ is $O_p(1)$. For more detailed descriptions of the processes $g_{c,d}$ and $g^0_{c,d}$, see [1, 3, 7, 26]. Thus, $g_{1,1} \equiv g$ and $g^0_{1,1} \equiv g^0$ are the unconstrained and constrained slope processes associated with the canonical process $X_{1,1}(z)$. Finally, define $\mathbb{D} := \int ((g(z))^2 - (g^0(z))^2)\,dz$.

The following theorem describes the limiting behavior of the unconstrained and constrained MLEs of $\psi$, appropriately normalized.

THEOREM 2.1. *Let*
$$X_n(h) = n^{1/3}(\hat{\psi}_n(z_0 + hn^{-1/3}) - \psi(z_0))$$
*and*
$$Y_n(h) = n^{1/3}(\hat{\psi}^0_n(z_0 + hn^{-1/3}) - \psi(z_0)).$$
*Let $a = (I(\psi(z_0))p_Z(z_0))^{-1/2}$ and $b = (1/2)\psi'(z_0)$. Under assumptions (A.0)–(A.7) and (a), (b), $(X_n(h), Y_n(h)) \to_d (g_{a,b}(h), g^0_{a,b}(h))$ finite dimensionally and also in the space $\mathcal{L} \times \mathcal{L}$.*

Thus $X_n(0) = n^{1/3}(\hat{\psi}_n(z_0) - \psi(z_0)) \to_d g_{a,b}(0)$. Using Brownian scaling it follows that the following distributional equality holds in the space $\mathcal{L} \times \mathcal{L}$:

(2.2) $\quad (g_{a,b}(h), g^0_{a,b}(h)) =_d (a(b/a)^{1/3} g((b/a)^{2/3} h), a(b/a)^{1/3} g^0((b/a)^{2/3} h)).$

For a proof of this proposition, see, for example, [1]. Using the fact that $g(0) \equiv_d 2\mathbb{Z}$ (see, e.g., [21]), we get

(2.3) $\qquad n^{1/3}(\hat{\psi}_n(z_0) - \psi(z_0)) \to_d a(b/a)^{1/3} g(0) \equiv_d (8a^2 b)^{1/3} \mathbb{Z}.$

This is precisely the phenomenon described in (1.1).

Our next theorem concerns the limit distribution of the likelihood ratio statistic for testing $H_0$.

THEOREM 2.2. *Under assumptions (A.0)–(A.7) and (a), (b),*
$$2 \log \lambda_n \to_d \mathbb{D} \quad \text{when } H_0 \text{ is true.}$$

REMARK 1. In this paper we work under the assumption that $Z$ has a Lebesgue density on its support. However, Theorems 2.1 and 2.2, the main results of this paper, continue to hold under the assumption that the distribution function of $Z$ is continuously differentiable (and hence has a Lebesgue density) in a neighborhood of $z_0$ with a nonvanishing derivative at $z_0$. Also, subsequently we tacitly assume that MLEs always exist; this is not really a stronger assumption than (A.0). Since our main results deal with convergence in distribution, we can, without loss of generality, restrict ourselves to sets with probability tending to 1. In this paper, we focus on the case where $\psi$ is increasing. The case where $\psi$ is decreasing is incorporated into this framework by replacing $Z$ by $-Z$ and considering the (increasing) function $\overline{\psi}(z) = \psi(-z)$.

*Characterizing $\hat{\psi}_n$*. In what follows, we define $\phi(x,\theta) \equiv -l(x,\theta)$, $\dot{\phi}(x,\theta) = -\dot{l}(x,\theta)$, $\ddot{\phi}(x,\theta) = -\ddot{l}(x,\theta)$ and $\phi'''(x,\theta) = -l'''(x,\theta)$. The log-likelihood function for the data is given by $\sum_{i=1}^{n} l(X_i, \psi(Z_i)) = \sum_{i=1}^{n} l(X_{(i)}, \psi(Z_{(i)}))$, where $Z_{(i)}$ is the $i$th smallest covariate value and $X_{(i)}$ is the response value corresponding to it. Finding the MLE under the constraint that $\psi$ is increasing reduces to minimizing $\tilde{\phi}(u_1, u_2, \ldots, u_n) = \sum_{i=1}^{n} \phi(X_{(i)}, u_i)$ over all $u_1 \leq u_2 \leq \cdots \leq u_n$. Once we obtain the (unique) minimizer $\hat{u} \equiv (\hat{u}_1, \hat{u}_2, \ldots, \hat{u}_n)$, the MLE $\hat{\psi}_n$ at the points $\{Z_{(i)}\}_{i=1}^{n}$ is given by $\hat{\psi}_n(Z_{(i)}) = \hat{u}_i$ for $i = 1, 2, \ldots, n$.

For convenience, take $\Theta$ to be $\mathbb{R}$ for the subsequent discussion (this assumption can be easily relaxed; see, in particular, Remark 2 below). By our assumptions, $\tilde{\phi}$ is a (continuous) convex function defined on $\mathbb{R}^n$ and necessary and sufficient conditions characterizing the minimizer are obtained readily, using the Kuhn–Tucker theorem. We write the constraints as $g(u) \leq 0$, where $g(u) = (g_1(u), g_2(u), \ldots, g_{n-1}(u))^T$ and $g_i(u) = u_i - u_{i+1}$, $i = 1, 2, \ldots, n-1$. Then there exists an $(n-1)$-dimensional vector $\lambda = (\lambda_1, \lambda_2, \ldots, \lambda_{n-1})^T$ with $\lambda_i \geq 0$ for all $i$, such that, if $\hat{u}$ is the minimizer satisfying the constraints, $g(\hat{u}) \leq 0$, then

$$\sum_{i=1}^{n-1} \lambda_i(\hat{u}_i - \hat{u}_{i+1}) = 0 \quad \text{and} \quad \nabla \tilde{\phi}(\hat{u}) + G^T \lambda = 0,$$

where $G$ is the $(n-1) \times n$ matrix of partial derivatives of $g$. The conditions displayed above are often referred to as *Fenchel conditions*. Solving recursively to obtain the $\lambda_i$'s (for $i = 1, 2, \ldots, n-1$), we get

$$\lambda_i \equiv \sum_{j=i+1}^{n} \nabla_j \tilde{\phi}(\hat{u}) = \sum_{i+1}^{n} \dot{\phi}(X_{(j)}, \hat{u}_j) \geq 0$$

(2.4)
$$\text{for } i = 1, 2, \ldots, (n-1)$$

and $\sum_{j=1}^{n} \nabla_j \tilde{\phi}(\hat{u}) = \sum_{j=1}^{n} \dot{\phi}(X_{(j)}, \hat{u}_j) = 0$. Now, let $B_1, B_2, \ldots, B_k$ be the blocks of indices on which the solution $\hat{u}$ is constant and let $w_j$ be the common value on block $B_j$. The equality $\sum_{i=1}^{n-1} \lambda_i(\hat{u}_i - \hat{u}_{i+1}) = 0$ forces $\lambda_i = 0$ whenever $\hat{u}_i < \hat{u}_{i+1}$. Noting that $\nabla_r \psi(\hat{u}) = \dot{\phi}(X_{(r)}, \hat{u}_r)$, this implies that on each $B_j$, $\sum_{r \in B_j} \dot{\phi}(X_{(r)}, w_j) = 0$. Thus $w_j$ is the unique solution to the equation $\sum_{r \in B_j} \dot{\phi}(X_{(r)}, w) = 0$. Also, if $S$ is a *head-subset* of the block $B_j$ (i.e., $S$ is the ordered subset of the first few indices of the ordered set $B_j$), then it follows that $\sum_{r \in S} \dot{\phi}(X_{(r)}, w_j) \leq 0$.

The solution $\hat{u}$ can be characterized as the vector of left derivatives of the greatest convex minorant (GCM) of a (random) cumulative sum (cusum) diagram, as will be shown below. The cusum diagram will itself be characterized in terms of the solution $\hat{u}$, giving us a *self-induced characterization*. Before proceeding further, we introduce some notation. For points $\{(x_i, y_i)\}_{i=0}^{n}$ where $x_0 = y_0 = 0$ and $x_0 < x_1 < \cdots < x_n$, consider the left-continuous function $P(x)$ such that $P(x_i) = y_i$ and such that $P(x)$ is constant on $(x_{i-1}, x_i)$. We will denote the vector of slopes (left-derivatives) of the GCM of $P(x)$ computed at the points $(x_1, x_2, \ldots, x_n)$ by $\text{slogcm}\{(x_i, y_i)\}_{i=0}^n$. Define the function

$$\xi(u) = \sum_{i=1}^n [u_i - \hat{u}_i + \nabla_i \tilde{\phi}(\hat{u}) d_i^{-1}]^2 d_i$$

$$= \sum_{i=1}^n [u_i - (\hat{u}_i - \dot{\phi}(X_{(i)}, \hat{u}_i) d_i^{-1})]^2 d_i,$$

where $d_i = \nabla_{ii}\tilde{\phi}(\hat{u}) = \ddot{\phi}(X_{(i)}, \hat{u}_i) > 0$. The function $\xi$ is strictly convex and it is easy to see that $\hat{u}$ minimizes $\xi$ subject to the constraints that $u_1 \leq u_2 \leq \cdots \leq u_n$ and hence, is given by the isotonic regression of the function $g(i) = \hat{u}_i - \dot{\phi}(X_{(i)}, \hat{u}_i) d_i^{-1}$ on the ordered set $\{1, 2, \ldots, n\}$ with weight function $d_i$. It is well known that the solution $(\hat{u}_1, \hat{u}_2, \ldots, \hat{u}_n) = \text{slogcm}\{\sum_{j=1}^i d_j, \sum_{j=1}^i g(j) d_j\}_{i=0}^n$. See, for example, Theorem 1.2.1 of [22]. In terms of the function $\phi$ the solution can be written as

(2.5)
$$\{\hat{u}_i\}_{i=1}^n \equiv \left[ \text{slogcm}\left\{ \sum_{j=1}^i \ddot{\phi}(X_{(j)}, \hat{u}_j), \right.\right.$$
$$\left.\left. \sum_{j=1}^i (\hat{u}_j \ddot{\phi}(X_{(j)}, \hat{u}_j) - \dot{\phi}(X_{(j)}, \hat{u}_j)) \right\}_{i=0}^n \right].$$

Recall that $\hat{\psi}_n(Z_{(i)}) = \hat{u}_i$; for a $z$ that lies strictly between $Z_{(i)}$ and $Z_{(i+1)}$, we set $\hat{\psi}_n(z) = \hat{\psi}_n(Z_{(i)})$. The MLE $\hat{\psi}_n$ thus defined is a piecewise constant right-continuous function.

*Characterizing* $\hat{\psi}_n^0$. Let $m$ be the number of $Z_i$'s that are less than or equal to $z_0$. Finding $\hat{\psi}_n^0$ amounts to minimizing $\tilde{\phi}(u) = \sum_{i=1}^n \phi(X_{(i)}, u_i)$ over all $u_1 \leq u_2 \leq \cdots \leq u_m \leq \theta_0 \leq u_{m+1} \leq \cdots \leq u_n$. This can be reduced to solving two separate optimization problems. These are: (1) Minimize $\sum_{i=1}^m \phi(X_{(i)}, u_i)$ over $u_1 \leq u_2 \leq \cdots \leq u_m \leq \theta_0$ and (2) Minimize $\sum_{i=m+1}^n \phi(X_{(i)}, u_i)$ over $\theta_0 \leq u_{m+1} \leq u_{m+2} \leq \cdots \leq u_n$.

Consider (1) first. As in the unconstrained minimization problem one can write down the Kuhn–Tucker conditions characterizing the minimizer. It is then easy to see that the solution $(\hat{u}_1^0, \hat{u}_2^0, \ldots, \hat{u}_m^0)$ can be obtained through the following recipe: Minimize $\sum_{i=1}^m \phi(X_{(i)}, u_i)$ over $u_1 \leq u_2 \leq \cdots \leq u_m$ to get $(\tilde{u}_1, \tilde{u}_2, \ldots, \tilde{u}_m)$. Then $(\hat{u}_1^0, \hat{u}_2^0, \ldots, \hat{u}_m^0) = (\tilde{u}_1 \wedge \theta_0, \tilde{u}_2 \wedge \theta_0, \ldots, \tilde{u}_m \wedge \theta_0)$. The solution vector to (2), say $(\hat{u}_{m+1}^0, \hat{u}_{m+2}^0, \ldots, \hat{u}_n^0)$, is similarly given by $(\hat{u}^0_{m+1}, \hat{u}^0_{m+2}, \ldots, \hat{u}^0_n) = (\tilde{u}_{m+1} \vee \theta_0, \tilde{u}_{m+2} \vee \theta_0, \ldots, \tilde{u}_n \vee \theta_0)$ where $\{\tilde{u}_i\}_{i=m+1}^n = \arg\min_{u_{m+1} \leq u_{m+2} \leq \cdots \leq u_n} \sum_{i=m+1}^n \phi(X_{(i)}, u_i)$.

An important property of the constrained solution $\{\hat{u}^0_i\}_{i=1}^n$ is that on any block $B$ of indices where it is constant and not equal to $\theta_0$, the constant value, say $w^0_B$, is the unique solution to the equation

$$\sum_{i \in B} \dot{\phi}(X_{(i)}, w) = 0. \tag{2.6}$$

The constrained solution also has a self-induced characterization in terms of the slope of the greatest convex minorant of a cumulative sum diagram. This follows in the same way as for the unconstrained solution by using the Kuhn–Tucker theorem and formulating a quadratic optimization problem based on the Fenchel conditions arising from this theorem. We skip the details but give the self-consistent characterization: The constrained solution $\hat{u}^0$ minimizes $A(u_1, u_2, \ldots, u_n) = \sum_{i=1}^n [u_i - (\hat{u}^0_i - \nabla_i \tilde{\phi}(\hat{u}^0) d_i^{-1})]^2 d_i$ subject to the constraints that $u_1 \leq u_2 \leq \cdots \leq u_m \leq \theta_0 \leq u_{m+1} \leq \cdots \leq u_n$, where $d_i = \nabla_{ii} \tilde{\phi}(\hat{u}^0)$. It is not difficult to see that

$$\{\hat{u}^0_i\}_{i=1}^m \equiv \left[ \text{slogcm} \left\{ \sum_{j=1}^i \ddot{\phi}(X_{(j)}, \hat{u}^0_j), \sum_{j=1}^i (\hat{u}^0_j \ddot{\phi}(X_{(j)}, \hat{u}^0_j) - \dot{\phi}(X_{(j)}, \hat{u}^0_j)) \right\}_{i=0}^m \right] \wedge \theta_0 \tag{2.7}$$

and

$$\{\hat{u}^0_i\}_{i=m+1}^n \equiv \left[ \text{slogcm} \left\{ \sum_{j=m+1}^i \ddot{\phi}(X_{(j)}, \hat{u}^0_j), \sum_{j=m+1}^i (\hat{u}^0_j \ddot{\phi}(X_{(j)}, \hat{u}^0_j) - \dot{\phi}(X_{(j)}, \hat{u}^0_j)) \right\}_{i=m}^n \right] \vee \theta_0. \tag{2.8}$$

The constrained MLE $\hat{\psi}^0_n$ is the piecewise constant right-continuous function satisfying $\hat{\psi}^0_n(Z_{(i)}) = \hat{u}^0_i$ for $i = 1, 2, \ldots, n$, $\hat{\psi}^0_n(z_0) = \theta_0$ and having no jump points outside the set $\{Z_{(i)}\}_{i=1}^n \cup \{z_0\}$.

REMARK 2. The characterization of the estimators above does not take into consideration boundary constraints on $\psi$. However, in certain models, the very nature of the problem imposes natural boundary constraints; for example, the parameter space $\Theta$ for the parametric model may be naturally nonnegative [example (d) discussed above], in which case the constraint $0 \leq u_1$ needs to be enforced. Similarly, there can be situations where $u_n$ is constrained to lie below some natural bound. In such cases, Fenchel conditions may be derived in the usual fashion by applying the Kuhn–Tucker theorem and self-induced characterizations may be derived similarly as above. However, as the sample size $n$ grows, with probability increasing to 1, the Fenchel conditions characterizing the estimator in a neighborhood of $z_0$ will remain unaffected by these additional boundary constraints, since $\psi(z_0)$ is assumed to lie in the interior of the parameter space, and the asymptotic distributional results will remain unaffected.

For this paper, we will assume the (uniform) almost sure consistency of the MLEs $\hat{\psi}_n$ and $\hat{\psi}_n^0$ for $\psi$ in a closed neighborhood of $z_0$. For the purpose of deducing the limit distributions of the MLEs and the likelihood ratio statistic, the following lemma, which guarantees local consistency at an appropriate rate, is crucial.

LEMMA 2.1. *For any $M_0 > 0$,*

$$\max\left\{\sup_{h\in[-M_0,M_0]}|\hat{\psi}_n(z_0+hn^{-1/3})-\psi(z_0)|,\right.$$
$$\left.\sup_{h\in[-M_0,M_0]}|\hat{\psi}_n^0(z_0+hn^{-1/3})-\psi(z_0)|\right\}$$

*is $O_p(n^{-1/3})$.*

We next state a number of preparatory lemmas required in the proofs of Theorems 2.1 and 2.2. But before that we need to introduce further notation. Let $\mathbb{P}_n$ denote the empirical measure of the data that assigns mass $1/n$ to each observation $(X_i, Z_i)$. For a monotone function $\Lambda$ defined on $\tilde{I}$ and taking values in $\Theta$, define the following processes: $W_{n,\Lambda}(r) = \mathbb{P}_n[\dot{\phi}(X,\Lambda(Z))1(Z\leq r)]$, $G_{n,\Lambda}(r) = \mathbb{P}_n[\ddot{\phi}(X,\Lambda(Z))1(Z\leq r)]$ and $B_{n,\Lambda}(r) = \int_{-\infty}^{r}\Lambda(z)\,dG_{n,\Lambda}(z) - W_{n,\Lambda}(r)$. Also, define normalized processes $\tilde{B}_{n,\Lambda}(h)$ and $\tilde{G}_{n,\Lambda}(h)$ in the following manner:

$$\tilde{B}_{n,\Lambda}(h) = n^{2/3}[(B_{n,\Lambda}(z_0+hn^{-1/3})-B_{n,\Lambda}(z_0))$$
$$-\psi(z_0)(G_{n,\Lambda}(z_0+hn^{-1/3})-G_{n,\Lambda}(z_0))]$$
$$\times (I(\psi(z_0))p_Z(z_0))^{-1}$$

and

$$\tilde{G}_{n,\Lambda}(h) = n^{1/3}\frac{1}{I(\psi(z_0))p_Z(z_0)}(G_{n,\Lambda}(z_0+hn^{-1/3})-G_{n,\Lambda}(z_0)).$$

LEMMA 2.2. *The process $\tilde{B}_{n,\psi}(h) \to_d X_{a,b}(h)$ in the space $B_{\mathrm{loc}}(\mathbb{R})$, where $a$ and $b$ are as defined in Theorem 2.1.*

LEMMA 2.3. *For every $K > 0$, the following asymptotic equivalences hold:*

$$\sup_{h \in [-K,K]} |\tilde{B}_{n,\psi}(h) - \tilde{B}_{n,\hat{\psi}_n}(h)| \to_p 0$$

*and*

$$\sup_{h \in [-K,K]} |\tilde{B}_{n,\psi}(h) - \tilde{B}_{n,\hat{\psi}_n^0}(h)| \to_p 0.$$

LEMMA 2.4. *The processes $\tilde{G}_{n,\hat{\psi}_n}(h)$ and $\tilde{G}_{n,\hat{\psi}_n^0}(h)$ both converge uniformly (in probability) to the deterministic function $h$ on the compact interval $[-K, K]$, for every $K > 0$.*

The next lemma characterizes the set $D_n$ on which $\hat{\psi}_n$ and $\hat{\psi}_n^0$ vary.

LEMMA 2.5. *Let $D_n$ denote the interval around $z_0$ on which $\hat{\psi}_n$ and $\hat{\psi}_n^0$ differ. Given any $\varepsilon > 0$, we can find an $M > 0$, such that for all sufficiently large $n$,*

$$P(D_n \subset [z_0 - Mn^{-1/3}, z_0 + Mn^{-1/3}]) \geq 1 - \varepsilon.$$

LEMMA 2.6 ([21]). *Suppose that $\{W_{n\varepsilon}\}, \{W_n\}$ and $\{W_\varepsilon\}$ are three sets of random vectors such that:*

(i) $\lim_{\varepsilon \to 0} \limsup_{n \to \infty} P[W_{n\varepsilon} \neq W_n] = 0$,
(ii) $\lim_{\varepsilon \to 0} P[W_\varepsilon \neq W] = 0$ *and*
(iii) *for every $\varepsilon > 0$, $W_{n\varepsilon} \to_d W_\varepsilon$ as $n \to \infty$.*

*Then $W_n \to_d W$, as $n \to \infty$.*

PROOF OF THEOREM 2.1. The proof presented here relies on continuous-mapping arguments for *slopes of greatest convex minorant* estimators. From the self-induced characterization of $\hat{\psi}_n$ [see (2.5)], we have

$$\hat{\psi}_n(z) - \psi(z_0) = \mathrm{slogcm}((B_{n,\hat{\psi}_n} - \psi(z_0)G_{n,\hat{\psi}_n}) \circ G_{n,\hat{\psi}_n}^{-1})(G_{n,\hat{\psi}_n}(z)).$$

Let $h \equiv n^{1/3}(z - z_0)$ be the local variable and recall the normalized processes that were defined before the statement of Lemma 2.2. In terms of the local variable and the normalized processes, it is not difficult to see that

$$n^{1/3}(\hat{\psi}_n(z_0 + hn^{-1/3}) - \psi(z_0)) = \mathrm{slogcm}(\tilde{B}_{n,\hat{\psi}_n} \circ \tilde{G}_{n,\hat{\psi}_n}^{-1})(\tilde{G}_{n,\hat{\psi}_n}(h)).$$

Similarly, from the characterization of $\hat{\psi}_n^0$ [refer to (2.7) and (2.8)] and the definitions of the normalized processes it follows that

$$n^{1/3}(\hat{\psi}_n^0(z_0 + hn^{-1/3}) - \psi(z_0)) = \text{slogcm}^0(\tilde{B}_{n,\hat{\psi}_n^0} \circ \tilde{G}_{n,\hat{\psi}_n^0}^{-1})(\tilde{G}_{n,\hat{\psi}_n^0}(h)).$$

Thus,

$$(X_n(h), Y_n(h)) = \{\text{slogcm}(\tilde{B}_{n,\hat{\psi}_n} \circ \tilde{G}_{n,\hat{\psi}_n}^{-1})(\tilde{G}_{n,\hat{\psi}_n}(h)), \\ \text{slogcm}^0(\tilde{B}_{n,\hat{\psi}_n^0} \circ \tilde{G}_{n,\hat{\psi}_n^0}^{-1})(\tilde{G}_{n,\hat{\psi}_n^0}(h))\}. \quad (2.9)$$

By Lemma 2.3, the processes $\tilde{B}_{n,\hat{\psi}_n^0}(h) - \tilde{B}_{n,\psi}(h)$ and $\tilde{B}_{n,\hat{\psi}_n}(h) - \tilde{B}_{n,\psi}(h)$ converge in probability to 0 uniformly on every compact set. Furthermore, by Lemma 2.2, the process $\tilde{B}_{n,\psi}(h)$ converges to the process $X_{a,b}(h)$ in $B_{\text{loc}}(\mathbb{R})$. It follows that the processes

$$(\tilde{B}_{n,\hat{\psi}_n^0}(h), \tilde{B}_{n,\hat{\psi}_n^0}(h)) \to_d (X_{a,b}(h), X_{a,b}(h)),$$

in the space $B_{\text{loc}}(\mathbb{R}) \times B_{\text{loc}}(\mathbb{R})$ equipped with the product topology. Furthermore, by Lemma 2.4, the processes

$$(\tilde{G}_{n,\hat{\psi}_n}(h), \tilde{G}_{n,\hat{\psi}_n^0}(h)) \to_p (h, h).$$

The proof is now completed by invoking continuous mapping arguments for slopes of greatest convex minorant estimators: thus, the limit distributions of $X_n$ and $Y_n$ are obtained by replacing the processes on the right-hand side of (2.9) by their limits. The details of the arguments are available in Theorem 2.1 of [2]. It follows that for any $(h_1, h_2, \ldots, h_k)$,

$$\{X_n(h_i), Y_n(h_i)\}_{i=1}^k \to_d \{\text{slogcm}\, X_{a,b}(h_i), \text{slogcm}^0 X_{a,b}(h_i)\}_{i=1}^k \\ = \{g_{a,b}(h_i), g_{a,b}^0(h_i)\}_{i=1}^k.$$

The above finite-dimensional convergence, coupled with the monotonicity of the functions involved, allows us to conclude that $(X_n(h), Y_n(h)) \to_d (g_{a,b}(h), g_{a,b}^0(h))$ in $\mathcal{L} \times \mathcal{L}$ as well. Indeed, if a sequence $\{\psi_n, \phi_n\}$ of monotone functions converges pointwise to the monotone functions $\{\psi, \phi\}$, then $(\psi_n, \phi_n)$ also converges to $(\psi, \phi)$ in $\mathcal{L} \times \mathcal{L}$ (see the result of Corollary 3 following Theorem 3 of [14]). □

PROOF OF THEOREM 2.2. We have

$$2\log \lambda_n = -2\left[\sum_{i \in J_n} \phi(X_{(i)}, \hat{\psi}_n(Z_{(i)})) - \sum_{i \in J_n} \phi(X_{(i)}, \hat{\psi}_n^0(Z_{(i)}))\right] \equiv -2S_n,$$

say. Here $J_n$ is the set of indices for which $\hat{\psi}_n(Z_{(i)})$ and $\hat{\psi}_n^0(Z_{(i)})$ are different. By Taylor expansion about $\psi(z_0)$, we find that $S_n$ equals

$$\left[\sum_{i \in J_n} \dot{\phi}(X_{(i)}, \psi(z_0))(\hat{\psi}_n(Z_{(i)}) - \psi(z_0))\right.$$
$$\left. + \sum_{i \in J_n} \frac{\ddot{\phi}(X_{(i)}, \psi(z_0))}{2}(\hat{\psi}_n(Z_{(i)}) - \psi(z_0))^2\right]$$
$$- \left[\sum_{i \in J_n} \dot{\phi}(X_{(i)}, \psi(z_0))(\hat{\psi}_n^0(Z_{(i)}) - \psi(z_0))\right.$$
$$\left. + \sum_{i \in J_n} \frac{\ddot{\phi}(X_{(i)}, \psi(z_0))}{2}(\hat{\psi}_n^0(Z_{(i)}) - \psi(z_0))^2\right] + R_n,$$

with $R_n = R_{n,1} - R_{n,2}$, where $R_{n,1} = (1/6)\sum_{i \in J_n} \phi'''(X_{(i)}, \psi_{n,i}^\star)(\hat{\psi}_n(Z_{(i)}) - \psi(z_0))^3$ and $R_{n,2} = (1/6)\sum_{i \in J_n} \phi'''(X_{(i)}, \psi_{n,i}^{\star\star})(\hat{\psi}_n^0(Z_{(i)}) - \psi(z_0))^3$ for points $\psi_{n,i}^\star$ [lying between $\hat{\psi}_n(Z_{(i)})$ and $\psi(z_0)$] and $\psi_{n,i}^{\star\star}$ [lying between $\hat{\psi}_n^0(Z_{(i)})$ and $\psi(z_0)$]. Under our assumptions $R_n$ is $o_p(1)$, as will be established later. Thus, we can write $S_n = I_n + II_n + o_p(1)$, where $I_n \equiv I_{n,1} - I_{n,2}$, with

$$(2.10) \quad I_{n,1} - I_{n,2} = \sum_{i \in J_n} \dot{\phi}(X_{(i)}, \psi(z_0))(\hat{\psi}_n(Z_{(i)}) - \psi(z_0))$$
$$- \sum_{i \in J_n} \dot{\phi}(X_{(i)}, \psi(z_0))(\hat{\psi}_n^0(Z_{(i)}) - \psi(z_0))$$

and

$$II_n = \sum_{i \in J_n} \frac{\ddot{\phi}(X_{(i)}, \psi(z_0))}{2}(\hat{\psi}_n(Z_{(i)}) - \psi(z_0))^2$$
$$- \sum_{i \in J_n} \frac{\ddot{\phi}(X_{(i)}, \psi(z_0))}{2}(\hat{\psi}_n^0(Z_{(i)}) - \psi(z_0))^2.$$

Consider the term $I_{n,2}$. Now, $J_n$ can be written as the union of blocks of indices, say $B_1^0, B_2^0, \ldots, B_l^0$, such that the constrained solution $\hat{\psi}_n^0$ is constant on each of these blocks. Let $B$ denote a typical block and let $w_B^0$ denote the constant value of the constrained MLE on this block; thus $\hat{\psi}_n^0(Z_{(j)}) = w_B^0$ for each $j \in B$. For any block $B$ where $w_B^0 \neq \theta_0$ we can write $\sum_{j \in B} \dot{\phi}(X_{(j)}, \psi(z_0)) \times (w_B^0 - \psi(z_0))$ as $(w_B^0 - \psi(z_0))H$, where

$$H = \sum_{j \in B}[\dot{\phi}(X_{(j)}, w_B^0) + (\psi(z_0) - w_B^0)\ddot{\phi}(X_{(j)}, w_B^0)$$
$$+ \tfrac{1}{2}(\psi(z_0) - w_B^0)^2 \phi'''(X_{(j)}, w_B^{0,\star})],$$

for some point $w_B^{0,\star}$ between $w_B^0$ and $\psi(z_0)$. Using the fact that on each block $B$ where $w_B^0 \neq \theta_0$, we have $\sum_{j \in B} \dot{\phi}(X_{(j)}, w_B^0) = 0$ [from (2.6)], it follows that $(w_B^0 - \psi(z_0))H$ equals

$$-\sum_{j \in B}(\psi(z_0) - w_B^0)^2 \ddot{\phi}(X_{(j)}, w_B^0)$$

$$-\tfrac{1}{2}\sum_{j \in B}(\psi(z_0) - w_B^0)^3 \phi'''(X_{(j)}, w_B^{0,\star}).$$

We conclude that $I_{n,2}$ equals

$$-\sum_{i \in J_n} \ddot{\phi}(X_{(i)}, \hat{\psi}_n^0(Z_{(i)}))(\hat{\psi}_n^0(Z_{(i)}) - \psi(z_0))^2$$

$$+ \tfrac{1}{2}\sum_{i \in J_n} \phi'''(X_{(i)}, \hat{\psi}_n^{0,\star}(Z_{(i)}))(\hat{\psi}_n^0(Z_{(i)}) - \psi(z_0))^3,$$

where $\hat{\psi}_n^{0,\star}(Z_{(i)})$ is a point between $\hat{\psi}_n^0(Z_{(i)})$ and $\psi(z_0)$. The second term in the above display is shown to be $o_p(1)$ by the exact same reasoning as used for $R_{n,1}$ or $R_{n,2}$. Hence, $I_{n,2} = -\sum_{i \in J_n} \ddot{\phi}(X_{(i)}, \hat{\psi}_n^0(Z_{(i)}))(\hat{\psi}_n^0(Z_{(i)}) - \psi(z_0))^2 + o_p(1)$, which by a one-step Taylor expansion about $\psi(z_0)$ can be seen to be equal to $-\sum_{i \in J_n} \ddot{\phi}(X_{(i)}, \psi(z_0))(\hat{\psi}_n^0(Z_{(i)}) - \psi(z_0))^2$ up to a $o_p(1)$ term. Similarly $I_{n,1} = -\sum_{i \in J_n} \ddot{\phi}(X_{(i)}, \psi(z_0))(\hat{\psi}_n(Z_{(i)}) - \psi(z_0))^2 + o_p(1)$. Now, using the fact that $S_n = I_{n,1} - I_{n,2} + II_n + o_p(1)$ and using the representations for these terms derived above, we find that up to a $o_p(1)$ term, $S_n$ equals

$$-\tfrac{1}{2}\bigg\{\sum_{i \in J_n} \ddot{\phi}(X_{(i)}, \psi(z_0))(\hat{\psi}_n(Z_{(i)}) - \psi(z_0))^2$$

$$-\sum_{i \in J_n} \ddot{\phi}(X_{(i)}, \psi(z_0))(\hat{\psi}_n^0(Z_{(i)}) - \psi(z_0))^2\bigg\},$$

whence $2 \log \lambda_n$ is given by

$$\sum_{i \in J_n} \ddot{\phi}(X_{(i)}, \psi(z_0))(\hat{\psi}_n(Z_{(i)}) - \psi(z_0))^2$$

$$-\sum_{i \in J_n} \ddot{\phi}(X_{(i)}, \psi(z_0))(\hat{\psi}_n^0(Z_{(i)}) - \psi(z_0))^2 + o_p(1).$$

Letting $\xi_n(x, z)$ denote the random function

$$\ddot{\phi}(x, \psi(z_0))\{(n^{1/3}(\hat{\psi}_n(z) - \psi(z_0)))^2 - (n^{1/3}(\hat{\psi}_n^0(z) - \psi(z_0)))^2\}1(z \in D_n),$$

it is easily seen that

$$2 \log \lambda_n = n^{1/3}(\mathbb{P}_n - P)\xi_n(x, z) + n^{1/3}P\xi_n(x, z) + o_p(1).$$

The term $n^{1/3}(\mathbb{P}_n - P)\xi_n(x, z) \to_p 0$ by Lemma 2.7 below. It now remains to deal with the term $n^{1/3} P(\xi_n(x, z))$ and as we will see, it is this term that contributes to the likelihood ratio statistic in the limit. We can write $n^{1/3} P \xi_n(x, z)$ as

$$n^{1/3} \int_{D_n} E_{\psi(z)}(\ddot{\phi}(X, \psi(z_0)))\{(n^{1/3}(\hat{\psi}_n(z) - \psi(z_0)))^2 - (n^{1/3}(\hat{\psi}_n^0(z) - \psi(z_0)))^2\} p_Z(z)\, dz.$$

On changing to the local variable $h = n^{1/3}(z - z_0)$ and denoting $z_0 + h n^{-1/3}$ by $z_n(h)$, the above can be decomposed as $A_n + B_n$, where

$$A_n \equiv \int_{\tilde{D}_n} [E_{\psi(z_0)} \ddot{\phi}(X, \psi(z_0))](X_n^2(h) - Y_n^2(h)) p_Z(z_n(h))\, dh$$

and

$$B_n \equiv \int_{\tilde{D}_n} [E_{\psi(z_n(h))} \ddot{\phi}(X, \psi(z_0)) - E_{\psi(z_0)} \ddot{\phi}(X, \psi(z_0))] \times (X_n^2(h) - Y_n^2(h)) p_Z(z_n(h))\, dh,$$

where $\tilde{D}_n = n^{1/3}(D_n - z_0)$. The term $B_n$ converges to 0 in probability on using the facts that eventually, with arbitrarily high probability, $\tilde{D}_n$ is contained in an interval of the form $[-M, M]$ on which the processes $X_n$ and $Y_n$ are $O_p(1)$ and that for every $M > 0$, $\sup_{|h| \leq M} | E_{\psi(z_0 + h n^{-1/3})}(\ddot{\phi}(X, \psi(z_0))) - E_{\psi(z_0)}(\ddot{\phi}(X, \psi(z_0))) | \to 0$ by (A.6). Thus,

$$2 \log \lambda_n = I(\psi(z_0)) \int_{\tilde{D}_n} (X_n^2(h) - Y_n^2(h)) p_Z(z_0 + h n^{-1/3})\, dh + o_p(1)$$

$$= \frac{1}{a^2} \int_{\tilde{D}_n} (X_n^2(h) - Y_n^2(h))\, dh + o_p(1).$$

We now deduce the asymptotic distribution of the expression on the right-hand side of the above display, using Lemma 2.6. Set $W_n = a^{-2} \int_{\tilde{D}_n} (X_n^2(h) - Y_n^2(h))\, dh$ and $W = a^{-2} \int \{(g_{a,b}(h))^2 - (g_{a,b}^0(h))^2\}\, dh$. Using Lemma 2.5, for each $\varepsilon > 0$, we can find a compact set $M_\varepsilon$ of the form $[-K_\varepsilon, K_\varepsilon]$ such that eventually, $P[\tilde{D}_n \subset [-K_\varepsilon, K_\varepsilon]] > 1 - \varepsilon$ and $P[D_{a,b} \subset [-K_\varepsilon, K_\varepsilon]] > 1 - \varepsilon$. Here $D_{a,b}$ is the set on which the processes $g_{a,b}$ and $g_{a,b}^0$ vary. Now let $W_{n\varepsilon} = a^{-2} \int_{[-K_\varepsilon, K_\varepsilon]} (X_n^2(h) - Y_n^2(h))\, dh$ and $W_\varepsilon = \int_{[-K_\varepsilon, K_\varepsilon]} (1/a^2)((g_{a,b}(h))^2 - (g_{a,b}^0(h))^2)\, dh$. Since $[-K_\varepsilon, K_\varepsilon]$ contains $\tilde{D}_n$ with probability greater than $1 - \varepsilon$ eventually ($\tilde{D}_n$ is the left closed, right open interval over which the processes $X_n$ and $Y_n$ differ), we have $P[W_{n\varepsilon} \neq W_n] < \varepsilon$ eventually. Similarly $P[W_\varepsilon \neq W] < \varepsilon$. Also $W_{n\varepsilon} \to_d W_\varepsilon$ as $n \to \infty$, for every fixed $\varepsilon$. This is so because by Theorem 2.1 $(X_n(h), Y_n(h)) \to_d (g_{a,b}(h), g_{a,b}^0(h))$ as a process in $\mathcal{L} \times \mathcal{L}$ and $(f_1, f_2) \mapsto \int_{[-c,c]} (f_1^2(h) - f_2^2(h)) \, dh$ is a continuous real-valued function defined from $\mathcal{L} \times \mathcal{L}$ to the reals. Thus all conditions of Lemma 2.6 are satisfied, leading to the conclusion that $W_n \to_d W$. The fact that the limiting distribution is actually independent of the constants $a$ and $b$, thereby showing universality, falls out from Brownian scaling. Using (2.2) we obtain

$$W = \frac{1}{a^2} \int \{(g_{a,b}(h))^2 - (g_{a,b}^0(h))^2\} \, dh$$

$$\equiv_d \frac{1}{a^2} a^2 (b/a)^{2/3} \int \{(g((b/a)^{2/3}h))^2 - (g^0((b/a)^{2/3}h))^2\} \, dh$$

$$= \int \{(g(w))^2 - (g^0(w))^2\} \, dh,$$

on making the change of variable $w = (b/a)^{2/3} h$. It only remains to show that $R_n$ is $o_p(1)$ as stated earlier. We outline the proof for $R_{n,1}$; the proof for $R_{n,2}$ is similar. We can write

$$R_{n,1} = \tfrac{1}{6} \mathbb{P}_n [\phi'''(X, \hat{\psi}_n^\star(Z)) \{n^{1/3}(\hat{\psi}_n(Z) - \psi(z_0))\}^3 1(Z \in D_n)],$$

where $\hat{\psi}_n^\star(Z)$ is some point between $\hat{\psi}_n(Z)$ and $\psi(z_0)$. On using the facts that $D_n$ is eventually contained in a set of the form $[z_0 - Mn^{-1/3}, z_0 + Mn^{-1/3}]$ with arbitrarily high probability on which $\{n^{1/3}(\hat{\psi}_n(Z) - \psi(z_0))\}^3$ is $O_p(1)$ and (A.5), we conclude that eventually, with arbitrarily high probability,

$$|R_{n,1}| \leq \tilde{C}(\mathbb{P}_n - P)[B(X) 1(Z \in [z_0 - Mn^{-1/3}, z_0 + Mn^{-1/3}])]$$
$$+ \tilde{C} P[B(X) 1(Z \in [z_0 - Mn^{-1/3}, z_0 + Mn^{-1/3}])],$$

for some constant $\tilde{C}$. That the first term on the right-hand side goes to 0 in probability is a consequence of an extended Glivenko–Cantelli theorem (see, e.g., Proposition 2 or Theorem 3 of [25]), whereas the second term goes to 0 by direct computation. $\square$

LEMMA 2.7. *With $\xi_n(x, z)$ as defined in the proof of Theorem 2.2, we have $n^{1/3}(\mathbb{P}_n - P)\xi_n(x, z) \to_p 0$.*

The proof of this lemma uses standard arguments from empirical process theory and can be found in [2].

**3. Applications of the main results.** In this section, we discuss some interesting special cases of monotone response models.

Consider the case of a one-parameter full rank exponential family model, naturally parametrized. Thus, $p(x, \theta) = \exp[\theta T(x) - C(\theta)] h(x)$ where $\theta$ varies in an open interval $\Theta$. The function $C$ possesses derivatives of all orders. Suppose we have $Z \sim p_Z(\cdot)$ and $X|Z = z \sim p(x, \psi(z))$ where $\psi$ is increasing or decreasing in $z$. We are interested in making inference on $\psi(z_0)$, where $z_0$ is an interior point in the support of $Z$. If $p_Z$ and $\psi$ satisfy conditions (a) and (b) of Section 2, the likelihood ratio statistic for testing $\psi(z_0) = \theta_0$ converges to $\mathbb{D}$ under the null hypothesis, since conditions (A.0)–(A.7) are readily satisfied for exponential family models. Note that $l(x, \theta) = \theta T(x) - C(\theta) + \log h(x)$, so that $\dot{l}(x, \theta) = T(x) - C'(\theta)$ and $\ddot{l}(x, \theta) = -C''(\theta)$. Conditions (A.1)–(A.6) can be checked quite readily. We leave the details to the reader. To check condition (A.7), note that since $C'(\theta)$ and $\ddot{l}(x, \theta)$ are uniformly bounded for $\theta \in (\theta_0 - \varepsilon, \theta_0 + \varepsilon) \equiv \mathcal{N}$, by choosing $M$ sufficiently large, we can ensure that for some constant $\gamma$ and $\theta \in \mathcal{N}$, $H(\theta, M) \leq E_\theta[(2T(X)^2 + \gamma)1\{|T(X)| > M/2\}]$, which in turn is dominated by

$$\sup_{\theta \in \mathcal{N}} e^{-C(\theta)} \left[ \int 2(T(x)^2 + \gamma)(e^{(\theta_0+\varepsilon)T(x)} + e^{(\theta_0-\varepsilon)T(x)}) \right.$$
$$\left. \times 1(|T(x)| > M/2) h(x) \, d\mu(x) \right].$$

The expression above is not dependent on $\theta$ and hence serves as a bound for $\sup_{\theta \in \mathcal{N}} H(\theta, M)$. As $M$ goes to $\infty$ the above expression goes to 0; this is seen by an appeal to the DCT and the fact that $T^2(X) + \gamma$ is integrable at parameter values $\theta_0 - \varepsilon$ and $\theta_0 + \varepsilon$.

The nice structure of exponential family models actually leads to a simpler characterization of the MLEs $\hat{\psi}_n$ and $\hat{\psi}_n^0$. For each block $B$ of indices on which $\hat{\psi}_n(Z_{(i)})$ is constant with common value equal to, say $w$, it follows from the discussion of the Fenchel conditions that characterize $\hat{\psi}_n$ in Section 2 that $\sum_{i \in B}(T(X_{(i)}) - C'(w)) = 0$; hence $C'(w) = n_B^{-1} \sum_{i \in B} T(X_{(i)})$ where $n_B$ is the number of indices in the block $B$. Furthermore, from the Fenchel conditions, it follows that if $S$ is a *head-subset* of the block $B$, then $\sum_{i \in B}(T(X_{(i)}) - C'(w)) \geq 0$; that is, $C'(w) \leq n_S^{-1} \sum_{i \in S} T(X_{(i)})$, $n_S$ denoting the cardinality of $S$. As a direct consequence of the above, we deduce that the unconstrained MLE $\hat{\psi}_n$ can actually be written as $\{C'(\hat{\psi}_n(Z_{(i)}))\}_{i=1}^n = \text{slogcm}\{G_n(Z_{(i)}), V_n(Z_{(i)})\}_{i=0}^n$, where $G_n(z) = n^{-1} \sum_{i=1}^n 1(Z_i \leq z)$ and $V_n(z) = n^{-1} \sum_{i=1}^n T(X_i) 1(Z_i \leq z)$ and $G_n(Z_{(0)}) \equiv V_n(Z_{(0)}) = 0$. The MLE $\hat{\psi}_n^0$ is characterized in a similar fashion but as constrained slopes of the cumulative sum diagram formed by the points $\{G_n(Z_{(i)}), V_n(Z_{(i)})\}_{i=0}^n$. Thus, the MLEs have explicit characterizations for these models and their asymptotic distributions may also be obtained by direct methods. It is not difficult to check that examples (a), (b) and (c) discussed in the Introduction are special cases of the one-parameter full rank exponential family models discussed above. Theorems 2.1 and 2.2 therefore hold for these models, and MLEs have explicit characterizations and are easily computable. In particular, the result on the likelihood ratio statistic derived in [3] follows as a special case of our current results.

We now turn to example (d); here, the conditional distribution of the response given the covariate comes from a curved exponential family. We study this example in some detail. In this case, with $c = d = 1$, the Fenchel conditions for $\hat{\psi}_n$ take the following form (in empirical process notation):

$$m \int_{[a,t)} \frac{1}{\hat{\psi}_n(z)} d\mathbb{P}_n(x,z) + \int_{z \in [a,t)} x\, d\mathbb{P}_n(x,z) + (m-1) \int_{z \in [a,t)} \frac{1}{\hat{\psi}_n(z)^{2m-1}}$$
$$- m \int_{z \in [a,t)} x^2 \hat{\psi}_n(z)^{2m-1}\, d\mathbb{P}_n(x,z) \geq 0, \qquad t \in [a,b],$$

with equality if $t$ is a jump point for $\hat{\psi}_n$. For a general (real) $m$, the above equations and inequalities do not translate into an explicit characterization of the MLE as the solution to an isotonic regression problem. The mileage that we get in the full rank exponential family models is no longer available to us owing to the more complex analytical structure of the current model. The self-induced characterization nevertheless allows us to write down the MLE as a slope of greatest convex minorant estimator, and determine its asymptotic behavior, following the route described in this paper (we leave the verification of the regularity conditions on the parametric model to the reader). To simplify matters, take $m = 1$. In this case, the joint distribution of $(X, Z)$ is given by $g(x, z) = p(x|\psi(z))h(z)$ where $h$ is the Lebesgue density of $Z$, $p(x|y) = yf(xy - 1)$ and $f$ is the standard normal density. For our simulation study, we chose $[a, b] = [1, 2]$ and $\psi(z) = z$ and $Z$ to follow the uniform distribution on $[1, 2]$. Then, $g(x, z) = p(x|z) = zf(xz - 1)$, for $x \in \mathbb{R}, z \in [1, 2]$. Essentially, we are in the setting of a mixture model. Following the discussion of the self-induced characterization for $\hat{\psi}_n$ described in Section 2, we find that $\hat{\psi}_n$, for this example, is the slope of the convex minorant of the "self-induced" cumulative sum diagram,

$$\{\mathbb{P}_n[(x^2 + \hat{\psi}_n(z)^{-2})1(z \in [a,t))], \mathbb{P}_n[(2\hat{\psi}_n^{-1}(z) + x)1(z \in [a,t))] : t \in [a,b]\}.$$

The left panel of Figure 1 shows the true function $\psi(z) = z$ and the MLE $\hat{\psi}_n(z)$ for the specific sample of size $n = 2000$ generated from the above model. It can be seen that the estimator tracks the true function quite well apart from the endpoints where the "spiking problem" manifests itself. The right panel of Figure 2 shows the unconstrained MLE (solid line) and the constrained MLE (dashed line) computed under the (true) $H_0 : \psi(1.5) = 1.5$, in a neighborhood of the point 1.5, along with the true function $\psi(z)$ (the slant line). The estimators are seen to coincide outside of a small interval around 1.5. It is the difference in behavior of the estimators in this short interval that contributes to the likelihood ratio statistic. The unconstrained estimator was obtained by running the ICM (iterative convex minorant) algorithm with starting value $\psi_0$ set to be the constant function 1.5. It converged quite rapidly without resorting to the line search procedure (see [15] for an excellent description of the modified ICM algorithm that incorporates line search to guarantee convergence). The constrained MLE was computed by decomposing the likelihood maximization procedure into two parts and optimizing separately, as described in the characterization of $\hat{\psi}_n^0$ in Section 2. The left panel of Figure 2 shows the quantile–quantile plot of 3000 points from the distribution of the likelihood ratio statistic for $n = 1000$ versus 3000 points from (a fine discrete approximation to) the distribution of $\mathbb{D}$, along with the line $y = x$. The quantile–quantile plot is in very good agreement with the line $y = x$, in conformity with the theory presented in this paper.

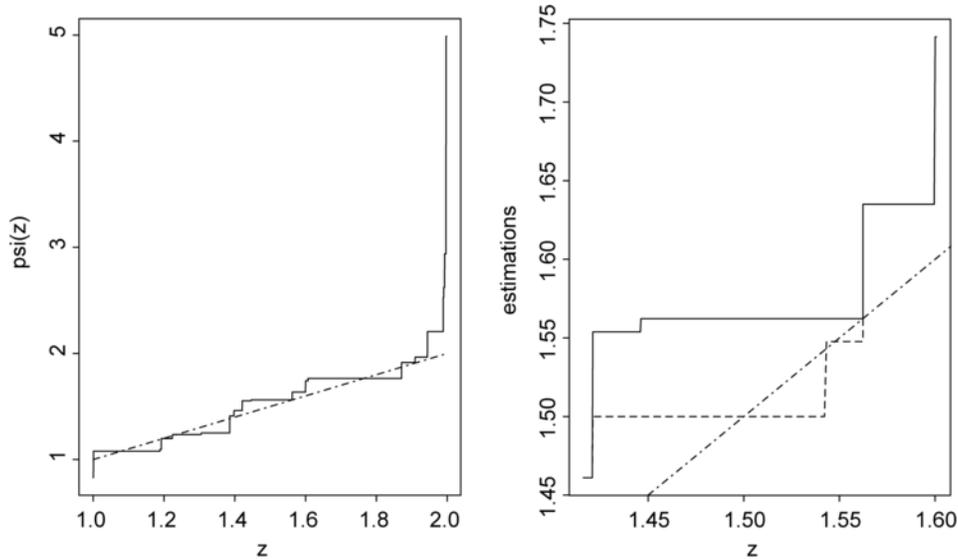

Fig. 1. *Left panel: The unconstrained estimator. Right panel: Close-up view of unconstrained and constrained estimators.*

An interesting fact that we now discuss is the *rapid* convergence of the ICM algorithm for this problem, in terms of number of steps to convergence. This is illustrated in the right panel of Figure 2. The histogram on the left is that of the number of iterations that is needed by the ICM to converge to $\hat{\psi}_n$ (with a tolerance of $10^{-5}$ for checking the Fenchel conditions) based on 1000 replicates for $n = 100$. The histogram on the right presents the same information but for $n = 10{,}000$. Despite the vast disparity in sample sizes, the histograms are very similar; less than 1% of the iterations consume more than ten steps (of course, the actual duration of an iteration is larger for $n = 10{,}000$). The fast convergence demonstrated through these histograms indicates that the performance of the ICM algorithm resembles that of the Newton algorithm, which is known to have good local convergence properties. This can be explained by the fact that the Hessian matrix of $\tilde{\phi}$ (minus the log-likelihood) is diagonal for this model, and indeed for the entire class of models considered in this paper, since the $u_i$'s, the arguments to $\tilde{\tilde{\phi}}$, are *separated* in the optimization problem (Indeed, the separation of variables that we encounter in the log-likelihood also allows us to solve the constrained optimization under $H_0$ by splitting the likelihood into two different parts and optimizing separately). A similar phenomenon was observed by Jongbloed [15] in his simulation studies on Case 2 interval censoring, where, despite the nondiagonal nature of the Hessian of the log-likelihood function (this is a *nonseparated* problem and will be discussed shortly), the ICM algorithm converged quickly, as a consequence of the fact that there were very few off-diagonal elements in the Hessian.

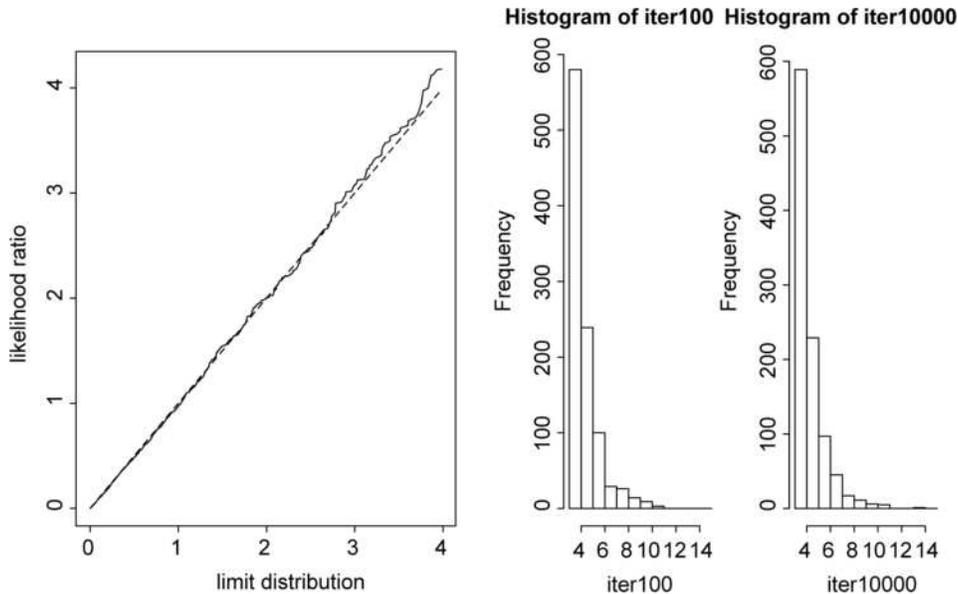

FIG. 2. *Left panel: Quantile–quantile plot of likelihood ratio statistic versus limiting quantiles. Right panel: Histograms of number of iterations until convergence.*

The Newton behavior of the ICM algorithm for the separated models of this paper suggests that in these models a one-step algorithm starting with the true function will produce estimators which are asymptotically equivalent to the MLE, even if the MLE is restricted by a null hypothesis. This phenomenon is alluded to as the "working hypothesis" in Section 5 of [9] and is illustrated through a derivation of the limit distribution of the MLE of the survival distribution $F$ for the current status model; in fact, the diagonal structure of the Hessian is used to establish the equivalence of the MLE with the "toy estimator" obtained by using the first iteration step of the ICM with the true distribution as the starting point. Thus, our results can be interpreted as pointing very strongly to the fact that the "working hypothesis" holds for the class of models considered in this paper.

While the approach in this paper applies nicely to separated models, there are several monotone function models of considerable interest where separatedness of the arguments to the log-likelihood function cannot be achieved; consequently the Hessian is no longer diagonal. Perhaps the simplest model of this type is the Case 2 interval censoring model, where there are two observation times $(U_i, V_i)$ for each individual, and one records in which of the three mutually disjoint intervals $(0, U_i]$, $(U_i, V_i]$, $(V_i, \infty)$ the individual fails. Letting $\{T_{(i)}\}_{i=1}^k$ denote the distinct ordered values of the $2n$ observation times $\{(U_i, V_i) : i = 1, 2, \ldots, n\}$ (here $n$ is the number of individuals being observed), and $u_l$ denote $F(T_{(l)})$, one can write down the log-likelihood for the data. It is seen that terms of the form $\log(u_i - u_j)$ immediately enter into the log-likelihood (see, e.g., [9], for a detailed treatment). One important consequence of this, in particular, is the fact that the computation of the constrained MLE of the survival distribution $F$ [under a hypothesis of the form $F(t_0) = \theta_0$] can no longer be decomposed into two separate optimization problems, in contrast to the monotone response models we have studied. Consequently, an analytical treatment of the constrained estimator will involve techniques beyond those presented in this paper. Regarding the unconstrained estimator of $F$ in this model, Groeneboom [8] uses some hard analysis to show that under a hypothesis of separation between $U$ and $V$ (the first and second observation times), the estimator converges to the truth at (pointwise) rate $n^{1/3}$, with limit distribution still given by $\mathbb{Z}$. The Case 2 model readily generalizes to the mixed case censoring model, where instead of two random observation times for every individual, the number of random times at which an individual is examined is also random. While heuristic considerations indicate that $\mathbb{D}$ should also arise as the limit distribution of the likelihood ratio statistic in these problems, the technical machinery for treating such nonseparated models in full generality remains to be developed and is left as a topic for future research.

**4. Discussion.** In this paper, we have studied the asymptotics of likelihood based inference in monotone response models. A crucial aspect of these models is the fact that conditional on the covariate $Z$, the response $X$ is generated from a parametric family that is regular in the usual sense; consequently, the conditional score functions, their derivatives and the conditional information play a key role in describing the asymptotic behavior of the maximum likelihood estimates of the function $\psi$. We have also shown that there are several monotone function models of interest that may be expected to exhibit asymptotically similar behavior though they are not monotone response models in the sense of this paper.

A potential extension of the monotone response models of this paper is to semiparametric models where the infinite-dimensional component is a monotone function. Here is a general formulation: Consider a random vector $(X, W, Z)$ where $Z$ is unidimensional, but $W$ can be vector-valued. Suppose that the distribution of $X$ conditional on $(W, Z) = (w, z)$ is given by $p(x, \beta^T w + \psi(z))$ where $p(x, \theta)$ is a one-dimensional parametric model. We are interested in making inference on both $\beta$ and $\psi$. The above formulation is fairly general and includes, for example, the partially linear regression model $X = \beta^T W + \psi(Z) + \varepsilon$ where $\psi$ is a monotone function (certain aspects of this model have been studied by Huang [12]), semiparametric logistic regression [with $X$ denoting the binary response and $(W, Z)$ covariates of interest] where the log odds of a positive outcome $(X = 1)$ is modeled as $\beta^T W + \psi(Z)$, and other models of interest. It is not difficult to see that the self-induced characterization will again come into play for describing the MLEs of $\psi$ in this general semiparametric setting. In the light of previous results, we expect that under appropriate conditions on $p(x, \theta)$, $\sqrt{n}(\hat{\beta}_{\mathrm{MLE}} - \beta)$ will converge to a normal distribution with asymptotic dispersion given by the inverse of the efficient information matrix and the likelihood ratio statistic for testing $\beta = \beta_0$ will be asymptotically $\chi^2$. The theory developed in [19, 20] should prove very useful in this regard. As far as estimation of the nonparametric component goes, $\hat{\psi}_n$, the MLE of $\psi$, should exhibit $n^{1/3}$ rate of convergence to a nonnormal limit and the likelihood ratio for testing $\psi$ pointwise should still converge to $\mathbb{D}$. This will be explored elsewhere and the ideas of the current paper should prove to be useful in dealing with the nonparametric component of the model.

## APPENDIX

Here, we present proofs of some selected lemmas. For proofs of the remaining lemmas, see [2].

PROOF OF LEMMA 2.2. It suffices to show that $\tilde{B}_{n,\psi}(h)$ converges to the process $aW(h) + bh^2$ in $l^\infty[-K, K]$, the space of uniformly bounded functions on $[-K, K]$ equipped with the topology of uniform convergence, for every $K > 0$. We can write

$$\tilde{B}_{n,\psi}(h) = \sqrt{n}(\mathbb{P}_n - P)f_{n,h} + \sqrt{n}Pf_{n,h},$$

where $f_{n,h}(X, Z)$ is given by

$$\frac{n^{1/6}[(\psi(Z) - \psi(z_0))\ddot{\phi}(X, \psi(Z)) - \dot{\phi}(X, \psi(Z))](1(Z \leq z_n(h)) - 1(Z \leq z_0))}{I(\psi(z_0))p_Z(z_0)},$$

with $z_n(h) \equiv z_0 + hn^{-1/3}$. To establish the above convergence, we invoke Theorem 2.11.22 of [24]. This requires verification of Conditions 2.11.21 and the convergence of the entropy integral in the statement of the theorem. Provided these conditions are satisfied, the sequence $\sqrt{n}(\mathbb{P}_n - P)f_{n,h}$ is asymptotically tight in $l^\infty[-K, K]$ and converges in distribution to a Gaussian process, the covariance kernel of which is given by
$$K(s,t) = \lim_{n \to \infty} (Pf_{n,s}f_{n,t} - Pf_{n,s}Pf_{n,t}).$$
We first compute $Pf_{n,s}f_{n,t}$. It is easy to see that this is 0 if $s$ and $t$ are of opposite signs, so we need only consider the cases where they both have the same sign. So let $s,t > 0$. Then, $Pf_{n,s}f_{n,t}$ is given by
$$E[n^{1/3}(\dot\phi(X,\psi(Z)) - (\psi(Z) - \psi(z_0))\ddot\phi(X,\psi(Z)))^2$$
$$\times 1(Z \in (z_0, z_0 + (s \wedge t)n^{-1/3}])] \times (I(\psi(z_0))p_Z(z_0))^{-2},$$
which can be written as $(I(\psi(z_0))p_Z(z_0))^{-2} n^{1/3} \int_{z_0}^{z_0+(s\wedge t)n^{-1/3}} G(\psi(z))p_Z(z)\,dz$, where, for every $\theta$, $G(\theta) = E_\theta[\dot\phi(X,\theta) - (\theta - \theta_0)\ddot\phi(X,\theta)]^2$. On expanding the square, $G(\theta)$ simplifies to
$$I(\theta) + (\theta - \theta_0)^2 f_3(\theta,\theta) - 2(\theta - \theta_0) E_\theta(\dot\phi(X,\theta)\ddot\phi(X,\theta)).$$
As $\theta \to \theta_0 \equiv \psi(z_0)$, the first term converges to $I(\theta_0)$ by (A.4) and the second term converges to 0 by (A.6). The third term also converges to 0, by the Cauchy–Schwarz inequality. It follows that $G(\theta)$ converges to $I(\theta_0) \equiv G(\theta_0)$. By the continuity of $\psi$ at $z_0$, we conclude that
$$\lim_{n\to\infty} Pf_{n,s}f_{n,t} = \frac{1}{(I(\psi(z_0)p_Z(z_0))^2} G(\psi(z_0))p_Z(z_0)(s \wedge t)$$
$$= \frac{1}{I(\psi(z_0))p_Z(z_0)}(s \wedge t).$$
It is easily shown that $Pf_{n,s}$ and $Pf_{n,t}$ both converge to 0 as $n \to \infty$, showing that for $s,t > 0$, $K(s,t) = [I(\psi(z_0))p_Z(z_0)]^{-1}(s \wedge t)$. Similarly, we can show that $K(s,t) = [I(\psi(z_0))p_Z(z_0)]^{-1}(|s| \wedge |t|)$, for $s,t < 0$. But this is the covariance kernel of the Gaussian process $aW(h)$ with $a = [I(\psi(z_0))p_Z(z_0)]^{-1/2}$. So the process $\sqrt{n}(\mathbb{P}_n - P)f_{n,h}$ converges in $l^\infty[-K, K]$ to the process $aW(h)$. We next show that $\sqrt{n}Pf_{n,h} \to (\psi'(z_0)/2)h^2$ uniformly on every $[-K, K]$. This implies that the process $\tilde{B}_{n,\psi}(h) \equiv \sqrt{n}\mathbb{P}_n f_{n,h}$ converges in distribution to $X_{a,b}(h) \equiv aW(h) + bh^2$ in $l^\infty[-K, K]$. To show the convergence of $\sqrt{n}Pf_{n,h}$ to the desired limit, we restrict ourselves to the case where $h > 0$; the case $h < 0$ can be handled similarly. Let $\xi_n(h) = I(\psi(z_0))p_Z(z_0)\sqrt{n}Pf_{n,h}$. Then $\xi_n(h)$ is given by
$$n^{2/3}E\{[(\psi(Z) - \psi(z_0))\ddot\phi(X,\psi(Z)) - \dot\phi(X,\psi(Z))]1(z_0 < Z \le z_0 + hn^{-1/3})\},$$
which reduces to $n^{2/3}E[(\psi(Z) - \psi(z_0))\ddot\phi(X,\psi(Z))1(z_0 < Z \le z_0 + hn^{-1/3})]$, on using the fact that $E_{\psi(z)}\dot\phi(X,\psi(z)) = 0$. Writing $z_n(u)$ for $z_0 + un^{-1/3}$ we can express this quantity as $A + B$ where $A = \int_0^h u\psi'(z_0)I(\psi(z_n(u)))p_Z(z_n(u))\,du$ and
$$B = \int_0^h [n^{1/3}(\psi(z_n(u)) - \psi(z_0)) - \psi'(z_0)u]I(\psi(z_n(u)))p_Z(z_n(u))\,du.$$

The term $B$ converges to 0 uniformly for $0 \leq h \leq K$ by the differentiability of $\psi$ at $z_0$ and $A$ can be written as $\int_0^h u\psi'(z_0)I(\psi(z_0))p_Z(z_0)\,du + o(1)$, where $o(1)$ goes to 0 uniformly over $h \in [0,K]$ and is readily seen to converge to $(1/2)(\psi'(z_0)I(\psi(z_0))p_Z(z_0))h^2$ uniformly on $0 \leq h \leq K$. It follows that $\sqrt{n}Pf_{n,h} \to (\psi'(z_0)/2)h^2$ uniformly over $0 \leq h \leq K$.

It remains to check Conditions 2.11.21. The computations here are tedious but straightforward, so we have omitted them (see [2] for the full details). Assumption (A.7), in particular, is used to verify a Lindeberg-type condition. $\square$

PROOF OF LEMMA 2.3. We only prove the first assertion. The second one follows similarly. For the first assertion, we write the proof for $h > 0$; the proof for $h < 0$ is similar. So, let $0 \leq h \leq K$. Recall that $\tilde{B}_{n,\hat{\psi}_n}(h)$ is given by

$$Cn^{2/3}\mathbb{P}_n[\{(\hat{\psi}_n(Z) - \psi(z_0))\ddot{\phi}(X, \hat{\psi}_n(Z)) - \dot{\phi}(X, \hat{\psi}_n(Z))\} \\ \times 1(Z \in (z_0, z_0 + hn^{-1/3}])],$$

where $C$ is a constant, and $\tilde{B}_{n,\psi}(h)$ has the same form as above but with $\hat{\psi}_n$ replaced by $\psi$. Now, for any $Z \in (z_0, z_0 + Kn^{-1/3}]$ we can write $\dot{\phi}(X, \psi(z_0))$ as

$$\dot{\phi}(X, \psi(Z)) + \ddot{\phi}(X, \psi(Z))(\psi(z_0) - \psi(Z)) \\ + \tfrac{1}{2}\phi'''(X, \psi^\star(Z))(\psi(Z) - \psi(z_0))^2,$$

for some point $\psi^\star(Z)$ between $\psi(Z)$ and $\psi(z_0)$. We can also write $\dot{\phi}(X, \psi(z_0))$ as

$$\dot{\phi}(X, \hat{\psi}_n(Z)) + \ddot{\phi}(X, \hat{\psi}_n(Z))(\psi(z_0) - \hat{\psi}_n(Z)) \\ + \tfrac{1}{2}\phi'''(X, \hat{\psi}_n^\star(Z))(\hat{\psi}_n(Z) - \psi(z_0))^2,$$

for some point $\hat{\psi}_n^\star(Z)$ between $\hat{\psi}_n(Z)$ and $\psi(z_0)$. It follows that we can write $\tilde{B}_{n,\psi}(h) - \tilde{B}_{n,\hat{\psi}_n}(h)$ as

$$C\tfrac{1}{2}\mathbb{P}_n[(n^{1/3}(\psi(Z) - \psi(z_0)))^2 \phi'''(X, \psi^\star(Z))1(Z \in (z_0, z_0 + hn^{-1/3}])] \\ - C\tfrac{1}{2}\mathbb{P}_n[(n^{1/3}(\hat{\psi}_n(Z) - \psi(z_0)))^2 \phi'''(X, \hat{\psi}_n^\star(Z))1(Z \in (z_0, z_0 + hn^{-1/3}])].$$

We will show that the second term in the above display converges to 0 uniformly in $h$; the proof for the first term is similar. Up to a constant, the second term is bounded in absolute value by

$$\mathbb{P}_n[(n^{1/3}(\hat{\psi}_n(Z) - \psi(z_0)))^2 |\phi'''(X, \hat{\psi}_n^\star(Z))|1(Z \in (z_0, z_0 + Kn^{-1/3}])].$$

Denote the random function inside square brackets by $\xi_n$. For any $z \in (z_0, z_0 + Kn^{-1/3}]$, we have

$$[n^{1/3}(\hat{\psi}_n(z) - \psi(z_0))]^2 \leq (n^{1/3}(\hat{\psi}_n(z_0 + Kn^{-1/3}) - \psi(z_0)))^2 \\ + (n^{1/3}(\hat{\psi}_n(z_0 - Kn^{-1/3}) - \psi(z_0)))^2,$$

which with arbitrarily high probability is eventually bounded by a constant $C$ (by Lemma 2.1). Also, since for any such $z$ $\hat{\psi}_n^\star(z)$ converges in probability to $\psi(z_0)$, with arbitrarily high probability $|\phi'''(X, \hat{\psi}_n^\star(Z))|$ is eventually bounded by $B(X)$ [by assumption (A.5)]. It follows that with arbitrarily high probability the random function $\xi_n$ is eventually bounded up to a constant by $B(X)1(Z \in [z_0, z_0 + Kn^{-1/3}])$. Hence, eventually, with arbitrarily high probability,

$$\mathbb{P}_n(\xi_n) \leq \tilde{C}(\mathbb{P}_n - P)[B(X)1(Z \in [z_0, z_0 + Kn^{-1/3}])] \\ + \tilde{C}P[B(X)1(Z \in [z_0, z_0 + Kn^{-1/3}])],$$

for some constant $\tilde{C}$. The first term on the right-hand side is $o_p(1)$ using straightforward Glivenko–Cantelli type arguments and the second term is seen to go to 0 by direct computation. This shows that the second term goes to 0 uniformly in $h$. □

**Acknowledgments.** I would like to thank Jon Wellner and Marloes Maathuis for their perceptive comments and some fruitful discussion. I would also like to acknowledge the inputs from the Associate Editor and two anonymous referees that resulted in a succinct but integrated presentation of the results in this paper.

Department of Statistics  
University of Michigan  
1085 South University  
Ann Arbor, Michigan 48109  
USA  
E-mail: moulib@umich.edu